\documentclass[11pt]{article}
\usepackage{graphicx}
\usepackage{amsmath}
\newcommand\lar[2]
  {\expandafter\gdef\csname labeled:#1\endcsname{#2}%
   \label{#1}#2}
\newcommand\recallLabel[1]
   {\csname labeled:#1\endcsname\tag{\ref{#1}}}

\usepackage{amssymb, amsmath}
\usepackage{epstopdf, wrapfig, framed,caption}
\usepackage{color}
\usepackage[makeroom]{cancel}
\usepackage{stmaryrd}
\usepackage{hyperref}
\usepackage{makecell}
\usepackage{authblk}

\usepackage{lineno}

\usepackage{newtxtext,newtxmath,amsmath}

\usepackage{etoolbox} 

\usepackage[makeroom]{cancel}

\def\bn{{\bf n}}
\def\bx{{\bf U}}
\def\be{{\bf e}}

\def\cL{{\cal L}}
\def\cP{{\cal P}}

\def\cR{{\mathcal R}}
\def\dz{{\dot{z}}}
\def\hS{{\hat{S}}}

\def\nV{\mathrm{V_n}}
\def\tnV{\tilde{\mathrm{V}}_{\mathrm n}}
\def\htheta{\hat{\theta}}

\def\hxi{\hat{\xi}}
\def\pn{\partial_\bn}
\def\pt{\partial_t}
\def\pz{\partial_z}
\def\salpha{}
\def\sg{\delta_g}
\def\R{\mathbb{R}}
\def\eps{\sigma_N}
\def\trd{\textrm d}
\def\opsi{\psi}
\def\oGam{\overline{\Gamma}}

\def\tn{\tilde{n}}
\def\tr{\tilde{\rho}}
\def\tt{\tilde{\theta}}
\def\tp{\tilde{\phi}}
\def\tP{\tilde{\Phi}}

\def\talpha{}
\def\tbeta{\beta}
\def\RR{\mathbb{R}}

\def\eN{{e_N}}
\def\eg{{e_g}}

\def\es{{e_s}}
\def\ewo{{e_1}}
\def\ewz{{e_0}}
\def\mrL{\mathrm L}
\numberwithin{equation}{section}

\def\beq{\begin{equation}}
\def\eeq{\end{equation}}

\newcommand{\norm}[1]{\left\lVert#1\right\rVert}

\newtoggle{HcScaling}
\toggletrue{HcScaling}

\title{Slow Migration of Brine Inclusions in First-Year Sea Ice}
\author[1]{Noa Kraitzman}
\author[2]{Keith Promislow}
\author[3]{Brian Wetton}
\affil[1]{Mathematical Sciences Institute, Australian National University, Australia, noa.kraitzman@anu.edu.au}
\affil[2]{Department of Mathematics, Michigan State University, United States, kpromisl@math.msu.edu}
\affil[3]{Department of Mathematics, University of British Columbia, Vancouver, Canada, wetton@math.ubc.ca}

\begin{document}

\maketitle
\begin{abstract}
  We derive a thermodynamically consistent model for phase change in sea ice by adding salt to the framework introduced by Penrose and Fife \cite{penrose1990thermodynamically}. Taking the salt entropy relative to the liquid water molar fraction provides a transparent mechanism for salt rejection under ice formation.  We identify slow varying coordinates, including salt density relative to liquid water molarity weighted by latent heat, and use multiscale analysis to derive a quasi-equilibrium Stefan-type problem via a sharp interface scaling. The singular limit is under-determined and the leading order system is closed by imposing local conservation of salt under interface perturbation. The quasi-steady system determines interface motion as  balance of curvature, temperature gradient, and salt density. We resolve this  numerically for axisymmetric surfaces and show that the thermal gradients typical of arctic sea ice can have a decisive impact on the mode of pinch-off of cylindrical brine inclusions and on the size distribution of the resultant spherical shapes. The density and distribution of inclusion sizes is a key component of sea ice albedo which factors into global climate models, \cite{perovich1996optical}.
\end{abstract}

{
  \small	
  \textbf{\text{Key words.}} Chemotaxis, Stefan problem, brine inclusions, thermodynamic self consistency,
}

 \section{Introduction}

Sea ice plays a significant role in ocean circulation and more broadly in the Earth's weather and climate system. At a global scale sea ice is a thin interfacial layer between the atmosphere and the ocean that serves to reduce evaporation, reflect sunlight, and insulate the ocean from heat loss. At a microscale sea ice is a highly complex, multi-component system composed of crystalline ice, liquid brine inclusions, air bubbles, and salt precipitate. This work focuses on brine inclusions. They play a leading role in sea ice microstructure, particularly determining its  electromagnetic and mechanical properties, creating habitat for a variety of CO$_2$-binding micro-algae \cite{junge2004bacterial}, and driving the oceanic transport of carbon, nutrients, and salt, \cite{golden2001brine}. Indeed, the fluid flow that arises from the onset of percolation of brine channels leads to both significant down-welling of brine laden water and an important increase in the effective thermal conductivity of sea ice \cite{Kraitzman2021Advection}. Sea ice's albedo, the percentage of solar radiation reflected, is greatly impacted by the surface area density of brine inclusions \cite{perovich1996optical}. Both albedo and effective thermal conductivity of sea ice are key parameters in global circulation models.

\begin{figure}
    \centering
    \begin{tabular}{cc}
    \includegraphics[width=3.5in]{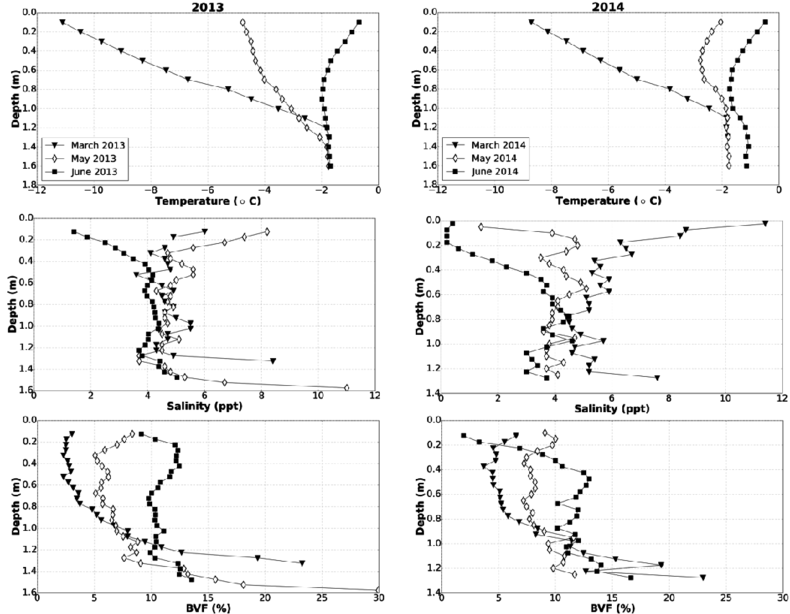}&
    \includegraphics[height=2in]{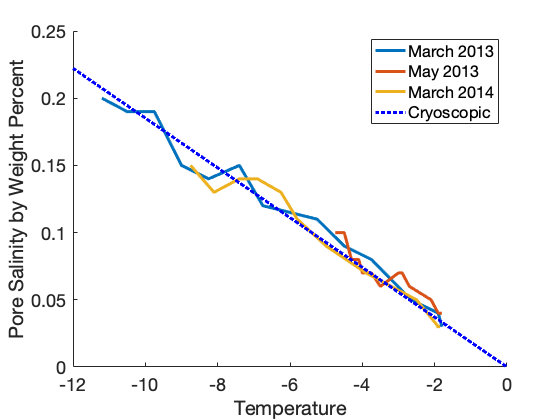}
    \end{tabular}
    \caption{(left) Measurements of temperature, bulk salinity (in parts per thousand), and brine volume fraction as a function of depth in arctic sea ice at three different seasonal periods over two years, \cite{o2016situ} reprinted with permission from \href{https://creativecommons.org/licenses/by/4.0/}{creative commons}. (right) Pore salinity computed as the ratio of bulk salinity and brine volume fraction for the three data sets (from left) with the largest temperature gradient. The pore salinity is plotted verses temperature (parametrically in depth) and compared to the linear approximation of the cryoscopic relation (dotted line).}
    \label{fig:Cryoscopic}
\end{figure}

Sea ice samples are widely characterized by their age -- first year or multiyear -- and the temperature and depth below the ice-air surface at which they were harvested.
The top row of images in Figure\,\ref{fig:Cryoscopic} (left), from \cite{o2016situ}, show the temperature distribution in sea ice at three different seasonal periods over two years. In a given season, the temperature is largely a function of depth while the temperature gradient is spatially uniform in the winter, except at the bottom (ice-water interface) of the ice sheet. The bulk salinity, shown in the second row, is the salt weight, in ppt, of the total volume, including the ice phase. The bulk salinity generically decays at the air-ice surface after the ice has been warmed, but remains relatively constant in the middle of the sheet. The brine volume fraction, calculated here from X-ray microtomography, increases with seasonal temperature shifts, but less so at the air-ice surface. While the data is noisy, ``pore salinity,'' the salt content of the brine within the inclusions, satisfies an underlying relation. The pore salinity can be post-processed from this data by approximating the salt density in the ice phase as zero, and dividing the bulk salt density by the brine volume fraction. The results, reported in  Figure\,\ref{fig:Cryoscopic} (right), present temperature and pore salinity parametrically by depth for the three data sets with the largest temperature gradient. Despite incorporating data from a wide range of brine structures throughout the ice-sheet, the plots are in strong agreement with the cryoscopic rule of thumb which relates a $0.54^\circ$C decrease in freezing temperature of water per $1$\% increase in salt weight (dotted line). The combination of salt ejection from ice and the raising of the freezing point of liquid water with salinity generates a chemotaxis process forming spatially extended liquid brine inclusions.

The structure of brine inclusions is well known to be temperature-dependent. As sea ice temperature rises, brine inclusions expand and may interconnect. As ice temperatures fall, inclusions shrink and pores may pinch-off into isolated inclusions of various shapes whose diameters range from $0.1$ to $10$ millimeters, \cite{light2003effects}, \cite{pringle2009pore}.
While there is an increasing interest in brine inclusions and sea ice structure  \cite{feltham2006sea, bartels2012ice, kumar2019generalized}, there are few mathematical models that incorporate the chemotactic role of salt ejection from ice and the cryoscopic relation between salt density and freezing point. There is a rich mathematical literature on phase change that incorporates latent heat, including the well cited \cite{caginalp1990dynamics} which presents scaling connections to many standard models. The work of \cite{fabrizio2016solidification} incorporates salt into a very general thermodynamic model with a strong emphasis on the elastic energy of the mixture, however this approach is technical and does not provide a transparent physical mechanism for salt ejection from ice.  The work \cite{morawetz2017formation} incorporates models for the nanoscale structure of water networks, and the role of salt in destabilizing them, but is not presented in the context of a thermodynamically consistent model.

We present a thermodynamically consistent phase-field model for the formation and evolution of brine inclusions within ice that rests upon a simple mechanism for salt exclusion.  The model follows the GENERIC framework for thermodynamic self consistency developed by Mielke \cite{mielke2011formulation}, incorporating salt into the entropy based models of phase change presented by Penrose and Fife, \cite{penrose1990thermodynamically}. Salt exclusion arises naturally by taking the entropy of the salt relative to the density of liquid-phase water molecules. The liquid water molecules solvate the salt ions, and their removal by the freezing process unfavorably decreases the entropy of the ions. The resulting ejection of salt from the regions of freezing engenders a chemotactic flow for the salt density that leads to the development of spatially extended regions of high salt concentration -- the brine inclusions.

We consider a scaling of the model which recovers a sharp-interface limit for the ice-liquid interface,  while maintaining a finite latent heat of phase change. The singular nature of the relative entropy of the salt precludes smooth transitions in salinity across the ice-water front. We overcome this by reformulating the system in terms of phase, temperature, and salt density relative to liquid phase, and show that this relative density is smooth.  We present a multiscale analysis that derives a Stefan-type problem (explicit moving boundary) for the evolution of the brine-ice front coupled to temperature, salinity, and interfacial curvature.
However because of the salt ejection the Stefan-type problem is formally underdetermined at leading order. We use conservation of salt under front perturbations to derive a boundary condition that closes the system. The Stefan-type problem is further simplified through a quasi-steady reduction to a mean curvature flow driven by the thermal gradient. This reduction aligns with recent work in the sea-ice community advocating for ``the removal of the widely adopted planar-equilibrium representation of the surface tension, the so-called capillary approximation, in favor of consideration of the curvature or size-dependence of the surface tension'', \cite{Hellmuth-2019}.

\begin{figure}
    \centering
    \includegraphics[width=2.5in]{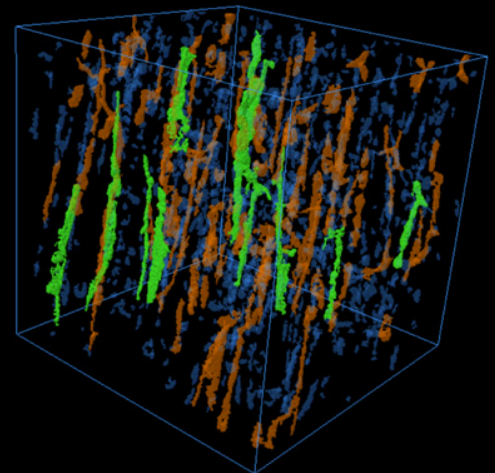}
    ~~
    \includegraphics[width=3.23in]{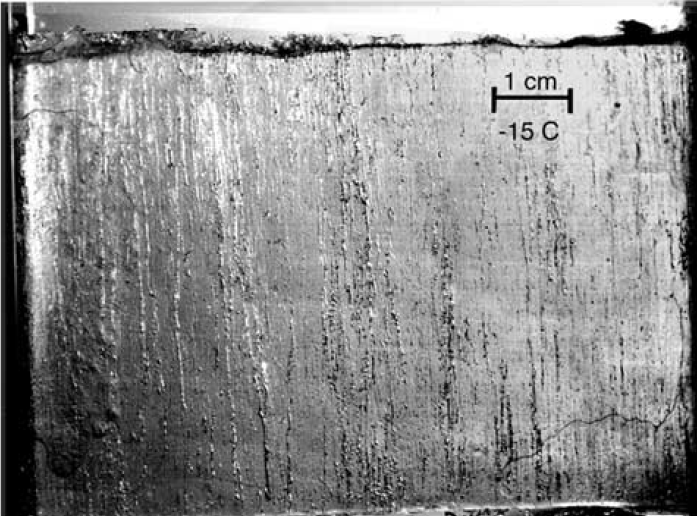}
    \caption{Vertical pores and smaller vertically aligned inclusions in first-year sea ice. (Left) Sea ice at -20$^\circ$C, scanned via X-ray microcomputed tomography. Edge of each side of cube is 7.5mm. Color coding indicates channel volume, green $> 0.0337$ mm$^3$, blue $< 0.0067$ mm$^3$, and orange are intermediate, \cite{lieb2017metrics}. (Right) A photograph of a vertical thick section cut from naturally grown first-year ice cut at a depth of approximately $80$ cm, \cite{light2003effects} reprinted with permission from John Wiley and Sons}.
    \label{fig:3D-brine}
\end{figure}

We calibrate the parameters to experimental data and use numerical simulations to examine the role of temperature gradients on pinch-off and migration of brine  inclusions.
We find that typical winter temperature gradients can have a significant impact on the structure of brine inclusions, contributing to the break-up (pinch-off) of longer brine pores into smaller, more spherical inclusions. More significantly for albedo, thermal gradients may impact the size distribution of the smaller inclusions. As shown in Figure \ref{fig:Wetton} a strong thermal gradient coupled with a large decrease in temperature may produce more uniformly sized inclusions after break-up, such as are visible in Figure\,\ref{fig:3D-brine}. The same thermal gradient with a smaller shift in temperature, as occurs towards the bottom of the ice sheet, produces a pinch-off at one end of the inclusion, suggesting that continuation beyond the initial pinch-off will lead to the formation of a string of small inclusions. The thermal gradient also induces a downward migration towards the warmer ice-bottom. This effect is particularly pronounced for spherical inclusions whose constant curvature cannot balance the inhomogeneity of the thermal gradient. This may provide a mechanism for the desalinization of the top layer of ice with the onset of warmer weather as is seen in the second row of data in Figure\,\ref{fig:Cryoscopic} (left). We investigate the contribution of the density stratification of salt within a pore to downward migration, but find that this effect is several orders of magnitude smaller.


There are important elements of sea ice that the model presented here does not capture.  It does not include air bubbles, it takes liquid water and ice to have the same density, ignoring expansion and the associated elastic effects induced by freezing. Sea ice is composed of a variety of salts, in particular mirabilite (Na$_2$SO$_4$) which begins to precipitate from solution at $-8.2^\circ$C and accumulate in the bottom of brine larger inclusions, impacting their shape and evolution, \cite{light2003effects}. Several effects,  including convective mixing during freezing, generate microstructure in sea ice that can influence brine channel shape, \cite{cole2001microstructure, lieb2017metrics}.  Core samples from the top of an ice sheet are more likely to be frazil ice with a random grain orientation, while mid-depth ice is generally columnar, with a more vertical orientation of ice grains, and ice near the seawater interface is often platelet ice, with a higher degree of lateral grain structure. Our model approach is best adapted to columnar ice.

In section 2 we derive the thermodynamically consistent flow, including the change of variables from temperature to internal energy, and the reformulation in terms of the more slowly varying salt density relative to liquid water. In section 3 we use multiscale analysis to pursue the sharp-interface limit, deriving the Stefan-type problem for the ice-water interface  and resolving the missing boundary condition. In section 4 we present the quasi-steady reduction to a thermal gradient driven curvature flow and examine the impact of the thermal gradient within the context of axisymmetric brine inclusions. Temperature is measured in $^\circ$K in sections 2 and 3, as required for discussions of thermal entropy which is singular at $0^\circ$K. In section 4 connections are made to experimental sea ice data and it is natural to change to $^\circ$C so that comparisons are more natural.

 \section{Model Derivation}

 We derive a thermodynamically consistent system for a mixture of water, ice, and salt within a cubical region $\Omega\subset{\mathbb R}^3$ subject  to zero-flux boundary conditions.
 This is a closed system  corresponding to a sample of ice at a fixed depth within an ice sheet. The evolution uniformly increases the system entropy, while preserving the internal energy and the total salt. In section\,\ref{s:Stefan} brine inclusion evolution in sea ice is modeled by replacing the zero-flux temperature boundary conditions on the top and bottom with inhomogeneous Dirichlet conditions that reflect the thermal gradient between the top surface that is cooled by the arctic air and the  bottom surface that is warmed by ocean water. With this driving force the system is no longer closed.

 The main variables are the phase function $\phi$, the temperature $\theta,$ and the salt weight fraction, $N$ (for Natrium). The phase function denotes the percentage of the water molecules that are in the  liquid state, thus $\phi=1$ is liquid, and $\phi=0$ is fully frozen.
 We incorporate salt to the phase field models introduced by Penrose and Fife in \cite{penrose1990thermodynamically} by following the GENERIC framework for thermodynamically consist models developed by Mielke, \cite{mielke2011formulation}. In the work of Penrose and Fife, the authors present several models for the free energy, internal energy, and entropy of temperature dependent phase change. We pursue the simplest of these, for which the entropy is expressed as a spatial integral of an entropy density in the form
\beq
s(\nabla\phi, \phi, \theta) = \nu(\phi) y(\theta,\phi)- \ewz W_0(\phi) - \sigma_\phi |\nabla\phi|^2.
\eeq
Here  $W_0$ is a classic double well potential with equal depth minima at $\phi=0,1$ and $\ewz$ and $\sigma_\phi$ are constants. While this formulation has some deficiencies, including a non-convex dependency of entropy on the phase field variable, they show that it connects directly to the classical phase field model and it significantly simplifies the subsequent analysis.
We modify their form, incorporating salt dependence into the first term in the form
$$
\nu(\phi)y(\theta,\phi, N) = - \ewo W_1(\phi)\xi(\theta,N).$$
The potential $W_1$ is dimensionless, and has $\phi=0,1$ as critical points so that
$$\tilde W(\phi,\theta, N):=\ewz W_0(\phi)+\ewo W_1(\phi)\xi(\theta,N),$$
is an unequal depth double well with respect to $\phi$ with local minima at $\phi=0,1$.
Generic choices are
\beq
\begin{aligned}
\label{e:W0-def}
W_0(\phi)&=18\phi^2(1-\phi)^2,\\
W_1(\phi)&=2 \phi^2
\left(\phi-\frac32\right).
\end{aligned}
\eeq
The cryoscopic term $\xi$ serves to raise or lower the value of the minima at $\phi=1$. For small deviations of temperature and low salt concentrations the cryoscopic term is well approximated by a linear relation,
 \beq
 \label{e:cryo-def}
 \xi(\theta, N):= \salpha N + \beta(\theta-\theta_*),
 \eeq
 where $\theta_*=273^\circ$K is the freezing point of pure water and $\beta=1.85/^\circ$K, see Figure\,\ref{fig:Cryoscopic} (right) and \cite{thompson1956concentration}. Since $W(0;\theta, N)=0$, the sign of $W(1;\theta, N)$ encodes the entropic preference of the mixture for ice or liquid. We take $W_1\leq 0$ on $[0,1]$, with the normalization $W_1(1)=-1$, so that positive values of $\xi$ promote melting. The normalization of $W_0$ simplifies the scaling of the surface tension, see \eqref{standingWave}.

The remaining impact of temperature and salt dependence is incorporated through the addition of two terms. The first is the simplest choice for the thermal entropy (see \cite{mielke2011formulation} page 238),
\beq \tilde\Gamma(\theta) = \es\ln\theta.
\eeq
The thermal entropy coefficient is the product of specific heat $c_s$ and density $\rho_0,$
$$\es=c_s\rho_0=(2500\,\textrm{J Kg}^{-1\circ} \textrm{K}^{-1} )(1000\,\textrm{Kg/m}^3)=2.5\times 10^6\, \textrm{J}^{\,\circ}\textrm{K}^{-1}\textrm{m}^{-3}.$$
The specific heat is taken independent of phase. The second term incorporates the entropy of the salt \emph{relative to the molar density of liquid} water. The water molecules in the liquid state solvate the salt ions. Consequently the salt entropy decreases with the ratio of salt molecules to liquid-state water molecules. This decrease in entropy drives the chemotactic ejection of salt from freezing water. Density driven stratification of salt within water is incorporated through a gravitational potential term,
 \beq
 \tilde\cR(N,\phi,x)= \eN  N\left(1-\ln\frac{N}{\phi}\right)-\eg Nx_3.
 \eeq
The entropy of NaCl salt in water (at $20^\circ$C) is roughly 43.4 J$^\circ$K$^{-1}$mol$^{-1}$ \cite{chemistry} while the molar density of water, $m_0= 5.55\times 10^4$ mol/m$^3$. This suggests an entropy coefficient $$\eN=(43.4\,\textrm{J}^\circ\textrm{K}^{-1}\,\textrm{mol}^{-1})\, m_0 =2.41\times 10^6\,\textrm{J}^\circ\textrm{K}^{-1}\textrm{m}^{-3}.$$
The coefficient $\eg$ quantifies the impact of gravity on the density of salt water relative to fresh water.
Water with $N$ percent salt by weight has density $\rho_0(1+0.83N)$ where $\rho_0$ is density of pure water, \cite{TheEngineeringToolBox}. This gives a buoyant density of the salt water of $\rho_0 0.83N$ with units of Kg/m$^3$. Introducing the gravitational constant $g$=9.8\,m/s$^2$ and scaling $x_3$ by the brine inclusion length $L_b=10^{-3}m$, we have
$$\eg := 0.83 \frac{g \rho_0 L_b}{\theta_*}=2.98 \times 10^{-2} \, \textrm{J}^\circ\textrm{K}^{-1} \textrm{m}^{-3}.$$
Here $\theta_*=273^\circ$K is the reference temperature for arctic salt water. The dimensionless ratio of salt entropy to its buoyant counterpart, akin to a Grashof number for thermal gradients, takes the value
$$\delta_g:=\frac{\eg}{\eN}= \frac{29.8\times 10^{-3}}{2.41\times 10^6} =1.23\times 10^{-8}.
$$
The constant $\ewo$ scales the latent heat. For water, the latent heat of freezing is $3.34\times 10^6$J/Kg, so that
$$\ewo=\frac{(3.34\times 10^6 \textrm{J Kg}^{-1}) \rho_0}{\theta_*} =1.2\times 10^7 \textrm{J}^\circ\textrm{K}^{-1}\textrm{m}^{-3}.$$
It is more difficult to estimate $\sigma_\phi$ and $\ewz$, however in the sharp interface regime the interfacial width satisfies $L_{li}=\sqrt{\ewz/\sigma_\phi}$. For a liquid-ice water interface this width is $L_{li}=10^{-9}$m. We define $H:=L_b/L_{li}\sim 10^6\gg 1$, which serves as the large parameter in our analysis. Non-dimensionalizing $x$ by the brine length-scale $L_b$ and the system entropy by $\es L_b^3$, we consider the classical sharp-interface scaling, imposing \beq
\begin{aligned}
\sigma_\phi&=\es L_{li}^2 H, \\
\ewz &=\es H.
\end{aligned}
\eeq
For simplicity of notation we set
$\eN=\ewo=\es.$
Dropping the tilde notation, the system entropy takes the form
 \beq
 S(\phi,\theta,N) =  \int_\Omega \overbrace{\Gamma(\theta) +\cR(N,\phi) - \frac{1}{H}|\nabla \phi|^2 - H W(\phi;\xi)}^{s(\nabla\phi,\phi,\theta,N)} \, \trd x.
 \eeq
The scaled potential takes the form
 \beq
 W(\phi;\theta, N) = W_0(\phi)+\frac{1}{H} W_1(\phi)\xi(\theta, N),
 \eeq
 while the thermal entropy
 \beq
 \Gamma=\ln \theta,
 \eeq
 and salt entropy relative to liquid water
 \beq
 \cR(N,\phi)= N\left(1-\ln \frac{N}{\phi}\right) -\delta_g N x_3,
 \eeq
 are dimensionless. We remark that time, measured in seconds, and temperature, measured in $^\circ$K, retain units.

 \subsection{The Entropic Gradient Flow}

We apply the thermodynamic framework of Mielke, see \cite{mielke2011formulation}, to develop a gradient flow that conserves the internal energy, increases the entropy, and conserves the total salt density. This requires replacing temperature with internal energy as a dependent variable. Assuming smoothness, we avoid the Legendre transform and express the free energy density $\opsi$  and internal energy density $u$ though the entropic density $s$,
 \beq s= -\partial_\theta \opsi,\eeq
 and
 \beq u = \opsi-\theta \partial_\theta \opsi = \opsi+\theta s.
 \eeq
 Up to terms that are independent of $\theta$ we compute that
 \beq\opsi = \frac{\theta}{2H} |\nabla \phi|^2 + \theta (H W_0(\phi) +\salpha NW_1(\phi)) + \frac{\beta(\theta-\theta_*)^2}{2} W_1(\phi) - \overline{\Gamma}(\theta) - \theta \cR(N,\phi).
 \eeq
 Here $\oGam$ is the primitive of $\Gamma$ with respect to $\theta$. The volume integral of the free energy density is a conserved quantity --  its gradients do not drive the flow. In this framework its significance is as an intermediate that determines the internal energy. This latter takes the form
 \beq
u =  -\frac{\beta}{2}(\theta^2-\theta_*^2)W_1(\phi) -\oGam(\theta)+\theta\Gamma(\theta),
\eeq
which we recast as
\beq
 \label{def-iE}
u = \left(\theta - b(\theta) W_1(\phi)\right),
\eeq
where we have introduced
\beq
\label{e:bdef}
b(\theta):= \frac{\beta}{2}(\theta^2-\theta_*^2).
\eeq
 Within this formulation the internal energy and entropy densities satisfy the fundamental thermodynamic relation
\beq
\label{TD-check}
\frac{\partial_\theta u}{\partial_\theta s} = \frac{1 -\beta\theta W_1(\phi)}{\theta^{-1}-\beta W_1(\phi)}= \theta.
\eeq

Subject to no-flux boundary conditions the thermodynamic evolution equations should increase the entropy density point-wise in space, while conserving the total internal energy
 \beq
 E_T:= \int_\Omega u(x)\,\trd x,
 \eeq
 and the total salt concentration
 \beq
 N_T:= \int_\Omega N(x)\,\trd x.
 \eeq
 This requires rewriting the evolution in terms of the internal energy and using  $(\phi,u, N)$ as dependent variables.  To this end we invert the relation \eqref{def-iE}, writing
 \beq
 \htheta:=\htheta(u,\phi).
 \eeq
 This inverse is well defined and smooth in the regime we consider.
The entropy takes the form
 \beq
 \label{def-hS}
 \hS(\phi, N, u) = \int_\Omega \Gamma(\htheta(u,\phi)) +\cR(N,\phi) - \frac{1}{2H}|\nabla \phi|^2 - H W(\phi;\hxi\,)\, \trd x,
 \eeq
 where to emphasize that the cryoscopic relation has become a function of $u$ and $\phi$  we introduce
 \beq
 \label{def-hxi}
 \hxi(\phi,u,N):=\salpha N+\beta(\htheta(\phi,u)-\theta_*).
 \eeq
The flow is determined by the gradient of $\hS$, through a choice of dissipation mechanism. The phase change is non-conservative, as ice and water interchange freely,
\beq
\phi_t = \frac{\delta \hS}{\delta \phi},
\eeq
where
\beq
\frac{\delta \hS}{\delta \phi} = \frac{1}{H}\Delta\phi- H W'(\phi,\hxi)+\left(\Gamma'(\htheta)-\beta W_1(\phi)\right)\partial_\phi\htheta+\partial_\phi \cR(N,\phi).
\eeq
Here prime denotes differentiation with respect to the function's dominant variable.
Taking $\partial_\phi$ of \eqref{def-iE} we calculate that
\beq
\partial_\phi \htheta  = \frac{\frac{\beta}{2}(\htheta^2-\theta_*^2)W_1'(\phi)}{1-\beta\htheta W_1(\phi)}.
\eeq
Since $\Gamma'=\theta^{-1}$ we find that the phase field equation for $\phi$ reduces to
\beq
\label{eq:phi1}
\phi_t= \frac{1}{H}\Delta\phi- HW'_0(\phi) -W'_1(\phi)\left(\salpha N+ B(\htheta) \right)+\partial_\phi \cR(N,\phi),
\eeq
where we have introduced $B(\theta)$, an increasing function of $\theta$ of the form
\beq\label{e-Bdef}
B(\theta):= \frac{\beta(\theta^2-\theta_*^2)}{2\theta}.
\eeq

For zero-flux boundary conditions the total internal energy is conserved by the flow
\beq
u_t = -\nabla \cdot\left( M_u \nabla\frac{\delta \hS}{\delta u}\right),
\eeq
where from \eqref{def-hS}
\begin{align}
\frac{\delta \hS}{\delta u} &= \left(\htheta^{-1}- \beta W_1(\phi)\right)\partial_u\htheta(u,\phi).
\end{align}
Taking $\partial_u$ of \eqref{def-iE} we derive
\beq
\partial_u\htheta = \frac{1}{1-\beta\htheta W_1(\phi)},
\eeq
and hence, as is consistent with \eqref{TD-check}, we obtain
\beq
\frac{\delta \hS}{\delta u} = \htheta^{-1}.
\eeq
For the canonical choice, $M_u=\sigma_\theta\htheta^2$, (see \cite{penrose1990thermodynamically} page 50), we derive the relation
\beq
\label{eq:u1}
u_t = \sigma_\theta \Delta\htheta,
\eeq
where $\sigma_\theta$ is the thermal conductivity, taken for simplicity to be independent of phase. The salt flux is generated by the entropy through the mass preserving flow
\beq
N_t = -\nabla \cdot\left( M_N \nabla \frac{\delta \hS}{\delta N} \right),
\eeq
where the variation of $\hS$ with respect to $N$ is given by
\begin{align}
\frac{\delta\hS}{\delta N} &=- W_1(\phi)\partial_N \hat\xi +\partial_NR(N,\phi)\\
& = -\salpha W_1(\phi)-\ln\left(\frac{N}{\phi}\right)- \delta_g{ x}_3.
\end{align}
Taking $M_N=\sigma_N N$, where $\sigma_N$ is a constant,  we obtain
\begin{align}
\label{eq:N1}
N_t &= \sigma_N\nabla\cdot \left( N \nabla \left(\salpha W_1(\phi)+\ln \frac{N}{\phi} +\delta_g{x}_3 \right) \right), \\
 &= \sigma_N\nabla \cdot\left(N\left(\salpha  W_1'(\phi) -\frac{1}{\phi}\right)\nabla\phi + \nabla N + \delta_g{\bf e}_3N \right).
\end{align}
To summarize, the coupled system takes the form
\begin{align}
    \phi_t &= \frac{1}{H}\Delta\phi- HW'_0(\phi) -W'_1(\phi)
    \left(\salpha N+B(\htheta)\right)
    +\frac{N}{\phi}, \label{e:GF-phi} \\
u_t  & = \sigma_\theta \Delta\htheta,
\label{e:GF-u}\\
N_t  &= \sigma_N\nabla \cdot\left(N\left(\salpha W_1'(\phi) - \frac{1}{\phi}\right)\nabla\phi + \nabla N+\delta_g {\bf e}_3 N \right).
\label{e:GF-N}
\end{align}

The advantage of this framework is that, subject to zero-flux boundary conditions, the system enjoys an entropy generation mechanism
\beq
 \frac{d\hS}{dt}=\int_\Omega \Big| \frac{\delta \hS}{\delta \phi}\Big|^2 +
 M_u \Big|\nabla \frac{\delta \hS}{\delta u} \Big|^2
 +
 M_N \Big|\nabla\frac{\delta \hS}{\delta N} \Big|^2 \,\trd x\geq 0,
 \eeq
 whose density is point-wise positive throughout the domain. In this sense the system satisfies the Clausius-Duhem (or second law) inequality, while maintaining conservation of total internal energy and salt. These conservation and entropy generation properties can be violated by the inhomogenous Dirichlet boundary conditions we impose on the the top and bottom of the domain in Section\,\ref{s:Stefan}. This does not violate the thermodynamic consistency as the system is no longer closed under inhomogeneous boundary conditions.

 \subsection{Relative Salt Density and Model Reformulation}
 For the multiscale analysis it is convenient to remove the fast variation from the salt variable, and to replace the internal energy with the temperature as a dependent variable.  The salt flux is given by
 \beq
 J_N:= N\sigma_N\left( \nabla \left(\salpha W_1+ \ln\frac{N}{\phi}\right)+\delta_g{\bf e}_3\right).
 \eeq
Setting aside the gravitational term, equilibrium with zero salt-flux have salt distributions of the form
 \beq
 N = \rho \phi e^{-\salpha W_1(\phi)},
 \eeq
 where $\rho>0$ is a spatial constant defining the salt density {\sl relative} to liquid water density.  This suggests that $\theta$ and $\rho$ can be viewed as slowly varying quantities, and we change to the variables $(\phi,\theta, \rho)$ from $(\phi, u, N)$. A key step is the introduction of the modified potential
 \beq\label{def-V}
 V(\phi;\theta,\rho) = W_0(\phi) +\frac{1}{H} V_1(\phi;\theta, \rho),
 \eeq
 where the  modified perturbation to the double well
 \beq
  \label{def-V1}
  \begin{aligned}
  V_1(\phi;\theta,\rho) &:= \int_0^\phi \left(B(\theta)+\salpha\rho s e^{-\salpha W_1(s)}\right) W_1'(s) -\rho e^{-\salpha W_1(s)}\, ds,\\
  &=  B(\theta)W_1(\phi) +\rho\left(\int_0^\phi -s\left(e^{-\salpha W_1(s)}\right)^\prime- e^{-\salpha W_1(s)}\,ds\right), \\
  & = B(\theta)W_1(\phi) -\rho\phi e^{-\salpha W_1(\phi)}.
  \end{aligned}
  \eeq
 This perturbation breaks the equal-depth structure, however, unlike $W_1$, it incorporates influence from the salt entropy, so that its value in the right well of $W_0$ is always negative. More specifically $V_1(0;\theta,\rho)=0$ while
 \beq\label{e:V1at1}
 V_1(1;\theta,\rho)=B(\theta)W_1(1) -\rho e^{-\salpha W_1(1)}=-B(\Theta)-\rho e^{1}<0.
 \eeq
 This shift in notation allows us to recast the system as
	\begin{align}
	\partial_t\phi &= \frac{1}{H}\Delta\phi-H \partial_\phi V(\phi;\theta,\rho),
	\label{e:SE2}\\
		\label{e:SE1}
	\partial_t\left(\theta -b(\theta) W_1(\phi)\right)&= \sigma_\theta \Delta \theta,\\
	\partial_t \left(\phi e^{-\salpha W_1(\phi)}\rho\right) & = \sigma_N \nabla\cdot \left(\phi e^{-\salpha W_1(\phi)} \left( \nabla \rho +\sg{\bf e}_3 \rho\right)\right). \label{e:SE3}
  \end{align}
The product of $W_1$ and the function $b(\theta)$, defined in \eqref{e:bdef}, provides the form and scaling for the latent heat.   In this formulation the salt exclusion mechanism is incorporated into the time derivative term on the left-hand side of \eqref{e:SE3}.

\begin{center}
Table of Parameters \\[2mm]
\begin{tabular}{ |c|l|c|c| }
\hline
Symbol & Name & Value & Units \\
\hline \hline
$\rho_0$ & Water density & 1000 & $\textrm{Kg/m}^3$\\
$c_s$ & Specific heat of ice & 2050 & $\textrm{J}/({^\circ}\hspace{-0.1em}\textrm{K}\,\textrm{Kg}) $\\
$S_e$ & NaCl salt in water molar entropy & 43.4 & $\textrm{J}/({^\circ}\hspace{-0.1em}\textrm{K}\,\textrm{mol})$\\
$m_0$ & Water molar density & $5.55\times 10^4$ & $\textrm{mol/m}^{3}$\\
$g$& Gravitational constant & 9.8 & $\textrm{m/s}^{2}$\\
$\theta_*$& Reference temperature & 273 & $^\circ\hspace{-0.1em}\textrm{K}$\\
$\beta$ & Cryoscopic parameter & $1.85$ & $1/^\circ$K\\
$e_s$ & Thermal entropy coefficient & $2.5\times10^6$ & $\textrm{J}/(^{\circ}\hspace{-0.1em}\textrm{K}\,\textrm{m}^{3})$\\
$e_N$ & Salt entropy coefficient & $2.41\times 10^6$  & $\textrm{J}/(^{\circ}\hspace{-0.1em}\textrm{K}\,\textrm{m}^{3})$\\
$e_g$ & Gravitational entropy coefficient & $2.98\times 10^{-2}$ & $\textrm{J}/(^{\circ}\hspace{-0.1em}\textrm{K}\,\textrm{m}^{3})$\\
$\ewo$ & Latent heat coefficient & $1.2\times 10^7$ & $\textrm{J}/(^\circ\textrm{K}\,\textrm{m}^{3})$\\
$L_b$  & Brine inclusion length scale & $10^{-3}$& $\textrm{m}$\\
$L_{li}$ & Liquid-ice interface length scale & $10^{-9}$ & $\textrm{m}$\\
$\delta_g$ & Density stratification ratio  & $1.23\times 10^{-8}$ & - \\
$H$ & Ratio of interface to inclusion lengths& $10^6$& - \\
 \hline
\end{tabular}
\end{center}

\section{Reduction to a Stefan-type Problem}

In the sharp-interface limit,  $H\gg 1$, we use multiscale asymptotics to derive a Stefan-type problem for the evolution of the ice-liquid interface. The salt rejection mechanism leads to non-smooth behavior in salt density across the interface, but the relative salt, $\rho$, remains smooth. The outer problem derived that results from matched asymptotics is under-determined, and we obtain a closed system by imposing conservation of salt under local interface variation.

We consider a brine inclusion whose boundary is given by a smooth, closed $2$-dimensional manifold $\Gamma$ immersed in $\Omega\subset\RR^3$. We define the local ``whiskered" coordinates system in a neighbourhood of $\Gamma$ via the mapping
\beq
    x = \zeta(p,z):=\gamma(p)+H^{-1}\bn(p)z,
\eeq
where $\gamma:\cP\rightarrow\RR^3$ is a local parameterization of $\Gamma$ and $\bn$ is the outward unit normal to $\Gamma$. The variables $p=(p_1,p_2)$ parameterize the tangential directions of $\Gamma$ while $z$ denotes the $H$-scaled, signed distance to $\Gamma$, negative on the region interior to $\Gamma$ and positive on the exterior. While we consider space dimension three for practical purposes, the arguments extend readily to co-dimension one inclusions in higher dimensions.

In a neighborhood $\Gamma_\ell$ of $\Gamma$
\beq
    \Gamma_\ell:=\left\{\zeta(p,z)\in\RR^3~\Big|~ p\in\cP ,~ -\ell H\leq z\leq \ell H\right\}\subset\Omega,
\eeq
the Cartesian gradient and Cartesian Laplacian admit the formulation
\begin{align}
    \nabla_x &= H\bn\partial_z+\nabla_p,\label{whiskers_grad}\\
    \Delta_x &= H^2\partial_z^2+ H \kappa_0(p)\partial_z +\left(z\kappa_1(p)\partial_z+\Delta_p\right)+O(H^{-1}).\label{whiskers_Lap}
\end{align}
Here $\kappa_i$ is the sum of the $i^{th}$ power of the curvatures
and, in particular, $\kappa_0$ is the total curvature of $\Gamma$. With this choice of normal the curvature of a sphere is negative, \cite{lee1997riemannian}.
The interface $\Gamma$ evolves in time through its normal velocity
 \beq\label{def:vel}
    \nV(p,t):=-H^{-1}\dz,
\eeq
 where $\dz:=\frac{\partial z}{\partial t}$ gives the rate of approach of the front to a point $x=x(z,p)$. It admits an expansion of the form
 \beq
    \dz = \dz_0+H^{-1}\dz_{-1}+O(H^{-2}).
 \eeq

 \subsection{Expansions and Matching Conditions}
The front profile $\tP$, is defined to be the solution of the equation
\beq\label{standingWave}
    \pz^2\tP = W_0'(\tP),
\eeq
which satisfies $\tP\rightarrow\tP^\pm$ as $z\rightarrow\pm\infty$ where $\tP^+=0$ and $\tP^-=1$ are the two minima of $W_0$. The linear operator obtained by linearizing \eqref{standingWave} about $\tP$ is denoted
\beq\label{def:L}
    \mrL:=\pz^2-W_0''(\tP).
\eeq
The function $\Phi\in H^2(\Omega)$ relates to $\tP$ through the relation $\Phi(x)=\tP(z(x))$ on $\Gamma_\ell$ and is extended smoothly to take values $\tP^\pm$ off of $\Gamma_\ell.$
We consider a formal, multiscale analysis of the temperature $\theta$, molar salinity $\rho,$ and the phase parameter $\phi$. In the far-field, away from the interface $\Gamma$, $\rho=N$ and the system admits the outer expansion
 \begin{align}
 \label{e:theta_exp}
    \theta &= \Theta_0(x,t)+H^{-1} \Theta_1(x,t)+O(H^{-2}),\\
    \phi &= \Phi (x,t)+H^{-1}\Phi_1(x,t)+O(H^{-2})
    \label{e:phi_exp},\\
    N &= N_0(x,t)+H^{-1} N_1(x,t)+O(H^{-2}).
    \label{e:N_exp}
 \end{align}
To avoid fast transients we consider a ``relaxed'' regime in which the phase change is equilibrated at leading order in the outer region. In the outer region, that is off of $\Gamma_\ell$, the leading order front profile $\Phi$ is piece-wise constant. This motivates the introduction of $\Omega_0(t)=\{x\, \bigl|\, \Phi (x,t)=0\}$ and $\Omega_1(t)=\{x\,\bigl| \,\Phi (x,t)=1\}$. We use $\chi_1$ the indicator of $\Omega_1$ to express $\Phi $, i.e.,
 $$
    \Phi (x,t) = \chi_1(x,t) = \begin{cases}
    1 & \text{if } x\in\Omega_1,\\
    0 & \text{if } x\in\Omega_0.
    \end{cases}
 $$
While the interface moves in time, this is reflected in the outer region only through the matching conditions.

 To simplify notation, we introduce the outer vector $\bx = [\Theta,N]^t$ and its expansion
\beq\label{e:vectorNotation}
\bx = \sum_{i\geq0} H^{-i}\bx_i,
\quad \bx_0 = [\Theta_0,N_0]^t,
\quad\bx_i=[\Theta_i,N_i]^t.
\eeq
We assume that the inner variables admit expansions of the form
 \begin{align}
    \theta(x,t) &=
    \tt(z,p,t) = \tt_0(z,p,t)+H^{-1} \tt_1(z,p,t)+O(H^{-2}),\label{innerTheta}\\
    \phi(x,t) &=
    \tp(z,p,t) = \tP(z)+ H^{-1}\tp_1(z,p,t)+O(H^{-2})
    \label{e:inner_phi_exp},\\
    \rho(x,t) &=
    \tr(z,p,t) = \tr_0(z,p,t)+H^{-1} \tr_1(z,p,t)+O(H^{-2}).\label{innerN}
 \end{align}
 At the interface we have matching conditions for both the temperature and the salinity.  The temperature satisfies the standard matching condition
 \begin{align}\label{match}
     \lim_{h\rightarrow\pm\infty}\Theta(x+h\bn,t) = \lim_{z\rightarrow\pm\infty}\tt(z,p,t),
 \end{align}
 which for $x\in\Gamma$ yields the relations
 \begin{align}
     \Theta_0^\pm(x,t) &= \lim_{z\rightarrow\pm\infty}\tt_0(z,p,t),\label{match_theta0}\\
     \Theta_1^\pm(x,t)+z\pn\Theta_0^\pm(x,t) & =
     \lim_{z\rightarrow\pm\infty}\tt_1(z,p,t).\label{match_theta1}
 \end{align}
 Here $\bn$ is the outward unit normal to $\Gamma$, $\pn$ is the derivative in the normal direction of $\Gamma$, and $\Theta_i^\pm$ denote the values of the limits of the left-hand side in \eqref{match} as $h\rightarrow\pm\infty$ respectively.

 The matching conditions on the salinity incorporate the relation $N = \rho \phi e^{-\salpha W_1(\phi)}$ and the usual limiting behavior,
  \begin{align}
     \lim_{h\rightarrow\pm\infty}N(x+h\bn,t) = \lim_{z\rightarrow\pm\infty}\tr(z,p,t) \tp(z,t,p) e^{-\salpha W_1(\tp(z,p,t))},
 \end{align}
 which for $x\in\Gamma$ yields the relations
 \begin{align}
     N_0^\pm(x,t) &= \lim_{z\rightarrow\pm\infty}\tr_0(z,p,t) \tP e^{-\salpha W_1(\tP)},\label{match_n0_TM}\\
     \label{match_r1_TM}
    N_1^\pm+z\pn N_0^\pm &= \lim_{z\to\pm\infty} \Big[\tP e^{-\salpha W_1(\tP)}(\tr_1-\tr_0\salpha W_1'(\tP)\tp_1)+\tp_1\tr_0e^{-\salpha W_1(\tP)}\Big].
\end{align}

\subsection{The Outer System}
We use the outer expansion \eqref{e:theta_exp}-\eqref{e:N_exp} to break the system \eqref{e:SE2}-\eqref{e:SE3} into orders of $H.$ In the ``relaxed'' outer regime, $\Phi=\chi_1$ for $x\in\Omega_0\cup\Omega_1$, and as a consequence the $O(H)$ system is trivially satisfied. The $O(1)$ equations take the form
\begin{align}
 W_0''(\chi_1)\Phi_1 &= -\partial_\phi V_1(\chi_1;\bx_0),\\
    \pt(\Theta_0-b(\Theta_0)W_1(\chi_1)) &= \sigma_\theta \Delta\Theta_0,\\
    \pt N_0  &= \sigma_N\nabla\cdot(\nabla N_0+\sg\be_3 N_0).
\end{align}
 Since $W_1(0) = 0$ and $W_1'(\chi_1)=0,$ the derivatives of the modified potential $V_1$ \eqref{def-V1} satisfy the relations
 \beq
     \partial_\phi V_1(0;\bx_0)=-\rho_0,\qquad
     \partial_\phi V_1(1;\bx_0)=-\rho_0e^{-\salpha W_1(1)}.
 \eeq
From the normalization $W_1(1)=-1$ we have $W_1(\chi_1)=-\chi_1$, and the $O(1)$ system can be written in terms of $\chi_1$ over the entire outer region as
\begin{align}
    \Phi_1 &= \frac{\rho_0e^{\chi_1}}{W''_0(\chi_1)}& \text{in } \Omega/\Gamma,\\
    \lar{TM:Theta_Outer}{(1+b(\Theta_0)\chi_1)\pt\Theta_0 &= \sigma_\theta\Delta\Theta_0,}\\
    \lar{TM:salt_Outer}{\pt N_0 &= \sigma_N\nabla\cdot(\nabla N_0+\sg\be_3N_0)},& \text{in } \Omega/\Gamma.
\end{align}

This system is subject to interior layer matching and exterior boundary conditions derived in the sequel. The phase parameter is zero in the ice region, and it is constant with a small perturbation which depends upon the salinity in the liquid region
\beq
    \phi(x,t) = 1+H^{-1}\frac{\rho_0 e^{-\salpha W_1(1)}}{W''_0(1)}+O(H^{-2}).
\eeq

\subsection{The Inner System}
In the inner region we combine the system \eqref{e:SE2}-\eqref{e:SE3}, the variable expansions \eqref{innerTheta}-\eqref{innerN}, and the gradient and the Laplacian expansions \eqref{whiskers_grad}-\eqref{whiskers_Lap}. Collecting  terms in orders of $H$ we find at $O(H^2)$
\begin{align}
    0&=\pz^2\tt_0,&\text{in }\Gamma_\ell\label{TM:tt_0}\\
    0&= \pz(\tP e^{-\salpha W_1(\tP)}\pz\tr_0),&\text{in }\Gamma_\ell.\label{TM:tr_0}
\end{align}
Equation \eqref{TM:tt_0} implies that $\tt_0$ is linear in $z$. The matching condition \eqref{match_theta0}  implies that $\tt_0$ is constant in $z$,  which yields the continuity condition on the outer temperature,
\beq\lar{TM:jump_t0}{
    \llbracket \Theta_0\rrbracket = 0,\quad \text{ on } \Gamma}.
\eeq
We conclude that $\tt_0=\Theta_0$. At $O(H)$, the system takes the form
\begin{align}
     0 &= \pz^2\tP-W_0'(\tP),\label{TM:tp_0}\\
     0 & = \sigma_\theta \pz^2\tt_1+\kappa_0\pz\tt_0,\label{TM:tt_1}\\
    \label{TM:tr_1}
    0 &= \pz(\tP e^{-\salpha W_1(\tP)}\pz\tr_1)+
    \left(\kappa_0\tP e^{-\salpha W_1(\tP)}+\pz(\tp_1e^{-\salpha W_1(\tP)}(1-\salpha W_1'(\tP)\tP))\right)\pz\tr_0+\\
    &\hspace{0.20in}\pz(\tP e^{-\salpha W_1(\tP)}) \sg\tr_0 \bn\cdot\be_3,\nonumber
\end{align}
%
The equation \eqref{TM:tp_0} and matching conditions are consistent with the assumption that $\tP$, provides the leading order inner expansion, in particular it is independent of the tangential variable $p$.

Since  $\tt_0$ is independent of $z$, equation \eqref{TM:tt_1} implies that $\tt_1$ is linear in $z$, and combined with the matching condition \eqref{match_theta1} yields the two interfacial zero-jump conditions for the outer temperature
\begin{align}
    \llbracket\pn\Theta_0\rrbracket =0,\quad
    \llbracket\Theta_1\rrbracket =
    0.
\end{align}

Addressing the $O(H^2)$ salt equation, \eqref{TM:tr_0}, we integrate twice with respect to $z$ from $0$ to $z$ and solve for $\tr_0$. This yields the relation
\beq
    \tr_0(z,p,t) =\tr_0(0)+ \left(\tP (0)e^{-\salpha W_1(\tP(0))}\partial_z\tr_0(0)\right)\int_0^z \frac{1}{\tP (y)e^{-\salpha W_1(\tP(y))}}\,dy.\label{TM:r0_sol}
\eeq
Reporting this back to the matching condition \eqref{match_n0_TM}, we have
\beq
N_0^\pm(x,t)= \lim_{z\to\pm\infty} \tP(z) e^{-\salpha W_1(\tP(z))}\left(
\tr_0(0)+ \left(\tP (0)e^{-\salpha W_1(\tP(0))}\partial_z\tr_0(0)\right)\int_0^z \frac{1}{\tP (y)e^{-\salpha W_1(\tP(y))}}\,dy
\right).
\eeq
  As $z\rightarrow -\infty$,  $\tP\rightarrow 1$ so that $W_1(\tP(z))\to W_1(1)$ and hence remains bounded.
 The dominant contribution comes from the term
\beq\label{TM:lim-} \lim_{z\rightarrow{-\infty}} e^{-\salpha W_1(\tP(z))} \tP(z)\int_0^z \frac{1}{e^{-\salpha W_1(\tP(y))}  \tP (y)}\,dy =
 z + O(1).
\eeq
The matching condition requires that $\pz\tr_0(0)= 0$, and hence $\tr_0$ is independent of $z$, and
\begin{align}
    &N_0^- = \tr_0(0)e^{-\salpha W_1(1)},&\text{on } \Gamma_\ell,
    \label{e:Match_N0-}\\
    &N_0^+ = 0,&\text{in } \Omega_0.
    \label{e:Match_N0+}
\end{align}
With $\tr_0$ independent of $z$, the $O(H)$ salt equation \eqref{TM:tr_1}
reduces to
\beq
    \pz(\tP e^{-\salpha W_1(\tP)}\pz\tr_1)
    =-\pz(\tP e^{-\salpha W_1(\tP)}) \sg\tr_0 \bn\cdot\be_3.
\eeq
This has solutions of the form
\beq
    \tr_1(z) = \tr_1(0)-z\sg\tr_0\bn\cdot\be_3+ \tP(0)e^{-\salpha W_1(\tP(0))}\left(\pz\tr_1(0)+\sg\tr_0\bn\cdot\be_3\right)\int_0^z\frac{1}{\tP(y)e^{-\salpha W_1(\tP(y))}}\,dy.
\eeq
We report this to the matching condition \eqref{match_r1_TM}.
Since $\tp_1$ is uniformly bounded and $\tr_0$ is independent of $z$ we use \eqref{TM:lim-}, to match terms in $z$ as $z\to-\infty$, finding that
\beq
\pn N_0^- = \tP(0)e^{-\salpha W_1(\tP(0))}\partial_z \tr_1 (0).
\eeq
Conversely, as $z\to\infty$ we have $\tP\to\tP^+= 0$, for which $W_1(0)=0$. Linearizing equation \eqref{TM:tp_0} about the limiting value $\tP^+$  yields the equation
 $$\pz^2\tP-W''(\tP^+)\tP = 0,$$
which implies that $\tP= c_+ e^{-k_+ z}$ as $z\to\infty$, where $k_+ = \sqrt{W''(\tP^+)}$. Using this asymptotic reduction, we have the relation
\beq\label{TM:lim+}
\begin{aligned}
 \lim_{z\to{+\infty}}  \tP(z,t)\int_0^z \frac{1}{\tP(y,t)}\,dy &=
\frac{1}{k_+}.
 \end{aligned}
\eeq
Using this limit in the matching condition we determine that
\beq
\pn N_0^+ = 0,
\eeq
which is consistent with \eqref{e:Match_N0+}.

The function $N_0$ is discontinuous across $\Gamma$, and  the interfacial conditions, expressed in terms of $\tr_0$ and $\partial_z\tr_1$, are under-determined. We close the interfacial condition for $N_0$ by imposing local conservation of salt mass under interface deformation. Since $N_0=0$ in the ice domain $\Omega_0$, conservation of mass requires that
$$\partial_t \int_{\Omega_1} N_0\, dx =0,$$
which, using \eqref{TM:salt_Outer} and the normal velocity $\nV$ breaks into
\begin{align}
0&= \int_\Gamma N_0 \nV\, ds + \int_{\Omega_1}\partial_t N_0\,dx,\\
&= \int_\Gamma N_0 \nV\, ds + \eps^2\int_{\Omega_1}\nabla\cdot(\nabla N_0+\sg\be_3 N_0)\,dx,\\
&= \int_\Gamma \left[N_0\nV+\eps^2\left( \nabla N_0+\sg N_0\be_3\right)\cdot\bn\right]\,ds.
\end{align}
Returning to the $\dot z$
formulation from \eqref{def:vel}, this implies the
leading order boundary condition
\begin{align}
\lar{TM:N_cons2}{&\frac{\dot{z}_0}{H} N_0 =\eps^2 \left(\nabla N_0+\sg N_0{\bf e}_3\right)\cdot \bn,& \text{on } \Gamma}.
\end{align}

A closed system for the leading order outer variables requires an expression for the normal velocity. This arises from the O(1) equation for the phase field variable which takes the form,
\beq\label{TM:tp1}
    \mrL\tp_1 = \dot{z}_0\tP' -\kappa_0\tP' - \partial_\phi V_1(\tP;\bx_0),
\eeq
where the operator $\mrL$ is defined in \eqref{def:L}. The solvability condition for \eqref{TM:tp1} requires that the right-hand side be orthogonal to the kernel $\tP'$ of $\mrL$. Taking the inner product and solving for $\dz_0$ yields the expression
\beq\label{TM:norm_vel_first}
    \dot{z}_0 = \kappa_0+\frac{\langle\partial_\phi V_1(\tP;\bx_0),\tP'\rangle}{\|\tP'\|_{L^2}^2}.
\eeq
Since $\bx_0$ is constant in $z$ it follows from \eqref{e:V1at1} that
$$
\left\langle\partial_\phi V_1(\tP;\bx_0),\tP'\right\rangle = V_1(\phi;\bx_0)\bigl|_{\phi=0}^{\phi=1} =
 B(\theta)W_1(1) -\rho e^{-\salpha W_1(1)}<0.$$

With the matching condition \eqref{e:Match_N0-} on $N_0^-$ and the normalization $W_1(1)=-1$, the expression  \eqref{TM:norm_vel_first} yields the scaled normal velocity
\begin{align}\lar{TM:norm_vel}{
    &\dot{z}_0 = \kappa_0-\norm{\Phi'}_{L^2}^{-2}\left( B(\Theta_0)+ N_0\right),
    & \text{on }\Gamma},
\end{align}
which couples temperature, salinity, and curvature. We summarize these results in section\,\ref{s:Stefan}.


\section{Stefan-type Problem for Brine Inclusions in Sea Ice}
\label{s:Stefan}
A brine inclusion is defined by its boundary $\Gamma$, which divides the scaled region $\Omega=[0,d_0]^3$ into subdomains $\Omega_0$ and $\Omega_1$. In the previous section we obtained a Stefan-type problem for the evolution of $\Gamma$ in terms of the leading order outer variables. The salt $N$ is zero on the exterior domain $\Omega_0$ and is discontinuous across the interface. It can be
taken to be defined only on the interior domain, $\Omega_1$. At leading order the phase equation is replaced by the location of the interface. We supplement the system with Dirichlet conditions on the temperature at the top $\partial \overline{\Omega}= [0,d_0]^2\times\{d_0\}$ and the bottom $\partial\underline{\Omega}:=  [0,d_0]^2\times\{0\}$, and
zero-flux conditions on the lateral sides $\partial\Omega_l:=\partial\Omega\backslash (\partial\overline{\Omega}\cup \partial\underline{\Omega})$.
The result is a nonlinear-parabolic equation for the temperature,
\begin{align}
   (1+b(\Theta_0)\chi_1)\pt\Theta_0 &= \sigma_\theta\Delta\Theta_0,& x\in \Omega\,,\\
  \Theta_0 &=\Theta_b(x,t), & \partial\overline{\Omega}\cup \partial\underline{\Omega},\,\\
   \partial_\bn \Theta_0 &=0, & x\in\partial\Omega_l,
  \end{align}
 where $\chi_1$ is the characteristic function for the inclusion region,  and a parabolic equation for the salt on the evolving brine inclusion region,
 \begin{align}
  %
  \pt N_0 &= \sigma_N\nabla\cdot(\nabla N_0+\sg\be_3N_0), & x\in \Omega_1,\\
         \partial_\bn N_0&=0, & x\in(\partial\Omega)\cap\Omega_1.
        \end{align}
The system is subject to the interfacial boundary condition which insures conservation of salt under the moving interface,
\begin{align}\label{TM:N_cons2_f}
\frac{\dot{z}_0}{H} N_0 &=\eps^2 \left(\nabla N_0+\sg N_0{\bf e}_3\right)\cdot \bn,& \text{on } \Gamma.
%
 \end{align}
 Here the $H$-scaled signed distance $z_0$ to $\Gamma$ satisfies
 \begin{align}
 &\dot{z}_0 = \kappa_0-\norm{\Phi'}_{L^2}^{-2}\left( B(\Theta_0)+ N_0\right),
    & \text{on }\Gamma.
 \label{e:normal-vel}
\end{align}
The system couples through the salt-preserving boundary condition, \eqref{TM:N_cons2_f} and the normal velocity \eqref{e:normal-vel}. The coupling is at $O(H^{-1})$, the same formal order as the second-order outer system. However on the long $O(H)$ time-scale both the leading and second order systems relax to quasi-equilibrium, and the coupling between the second order and first order reduces to $O(H^{-2}),$ and is negligible. The leading order temperature system is parabolic with nonlinearity arising only through the temperature and spatial dependence of the latent heat, $b(\theta)\chi_1$. Since $b>0$, the system is uniformly parabolic in non-divergence form, and the parabolic regularity theory, see Section 7.1 of \cite{evans1998partial}, applies.
The salt weight fraction satisfies a scaled advection-diffusion equation within $\Omega_1$ and is defined to be zero outside this domain. Standard regularity theory applies to this system too.

\subsection{Quasi-equilibrium Stefan-type problem and its axisymmetric formulation}

As can be seen from \eqref{TM:N_cons2_f} on the fast $t$ time scale the domain is constant to leading order, and
both the heat and the salt equations satisfy unforced parabolic equations. Thus these quantities relax to quasi-equilibrium on this time-scale, and are driven adiabatically by the interface which evolves on a slower $\tau=t/H$ time scale. In particular, for temperature boundary data that is spatially uniform on the top, $\partial\overline{\Omega}$, and the bottom, $\partial\underline{\Omega}$, the temperature relaxes to a simple linear equilibrium
\beq
\label{e:Theta0-steady}
\Theta_0=a_0+b_0x_3,
\eeq
with $b_0<0$ reflecting that sea ice is generically warmer with increasing depth. Assuming that the brine inclusion region $\Omega_1$ does not intersect $\partial\Omega$, the salt density satisfies
\begin{align}
\label{e:N0-steady}
N_0 &= N_T\frac{e^{-\sg x_3}}{\int_{\Omega_1} e^{-\sg x_3}\,dx}&x\in\Omega_1,
\end{align}
where the normalization incorporates the conservation of total salt, $N_T$.

The front evolution is driven quasi-adiabatically through the normal velocity. For simplicity we linearize $B(\Theta)$, defined in \eqref{e-Bdef}, about $\Theta_0=\theta_*$, consider the  slow time scale, $\tau=t/H$ and observe that with the normalization \eqref{e:W0-def} we have $\|\tP'\|_{L^2}=1.$ With these adjustments the inner normal velocity $\tnV=\nV/H$ takes the form
\beq
\label{e:S-NV}
\tnV=-\kappa_0+ \xi(N_0,\Theta_0)= -\kappa_0+ N_0+\beta (\Theta_0-\theta_*).
\eeq
The linearized normal velocity balances the cryoscopic term against curvature.
Since $H\sim 10^6,$ the slow time $\tau=1$ corresponds to roughly two weeks, and $\tau=5-10$ comprises a complete winter season.


 With $\Theta_0$ and $N_0$ prescribed as in \eqref{e:Theta0-steady}-\eqref{e:N0-steady}, brine inclusions achieve shape equilibrium when the curvature balances the heat and salt gradients. Given that brine inclusions are predominantly spherical and vertically oriented cylindrical pores, it is natural to consider a vertically oriented axisymmetric reduction for the curvature flow. An axisymmetric surface of revolution has a parameterization
\beq\label{e:axisym-param} \sigma(s,\mu)=(r(s)\cos\mu,r(s)\sin\mu, s),
\eeq
over $[0,d]\times[0,2\pi]$ where $r:[0,d]\mapsto\mathbb{R}$ denotes the radius of the surface measured from its vertical center line.  The curvature relates to $r$ through the equality
\beq\label{e:curvature} \kappa_0 =- \frac{1+(r')^2 -r r''}{2r\sqrt{(1+(r')^2)^3}}.
\eeq
For zero normal velocity the curvature satisfies
\beq
\kappa_0= \talpha N_T\frac{ e^{-\sg x_3}}{\int_a^b \pi r^2 e^{-\sg x_3}\,dx_3} +\tbeta(a_0+b_0x_3-\theta_*),
\eeq
which can be solved as a second order ODE for $r=r(x_3).$

More generally, under evolution by a normal velocity $\nV$ the the map $x_3=x_3(s)$ becomes non-trivial and the axisymmetric parameterization takes the form
$$ \sigma(s,\mu)=(r(s)\cos\mu,r(s)\sin\mu, x_3(s)).$$
The outer normal to the interface is given by
$$\bn = \frac{\left(x_3'\cos\mu,x_3'\sin\mu,-r'\right)}{\sqrt{(x_3')^2+(r')^2}},$$
while the time dependent $r$ and $x_3$ parameterizations satisfy
\beq
\begin{aligned}
\label{e:asym-NV}
\pt r &= \frac{\tnV x_3' }{\sqrt{(x_3')^2+(r')^2}}, \\
\pt x_3& = -\frac{\tnV r'}{\sqrt{(x_3')^2+(r')^2}}.
\end{aligned}
\eeq
The curvature satisfies
$$ \kappa_0(s) =- \frac{(x_3')^3+(r')^2x_3' -rr''x_3'+rr'x_{3}''} {2r\sqrt{((x_3')^2+(r')^2)^3}},$$
where $\prime$ denotes $\partial_s.$ The normal velocity is then computed in terms of $s$ as
\beq
\label{e:NV-s}
\tnV(s)= -\kappa_0(s)+\talpha N_0(x_3(s))+\tbeta(\Theta_0(x_3(s))-\theta_*).
\eeq
Returning \eqref{e:NV-s} to \eqref{e:asym-NV} gives a closed evolution for $r(s,t)$ and $x_3(s,t)$ on the fixed domain $[0,d].$

This axially symmetric dynamic problem is computed on a cell-centered grid with $s \in [0,1]$ a scaled arc length variable as in \cite{discretization}. Finite difference approximations are used for the derivatives and the integral in \ref{e:N0-steady} is approximated with the trapezoidal rule. The right half of the shape is computed and ghost points \cite{ghost} are used to apply the boundary and symmetry conditions. Backward Euler time stepping is implemented and Newton iterations are performed at every time step to solve the resulting nonlinear system. A grid refinement study gave the expected convergence: first order accuracy in time-step and second order accuracy with respect to spatial resolution.



\begin{figure}
    \centering
        \includegraphics[width =6.5in]{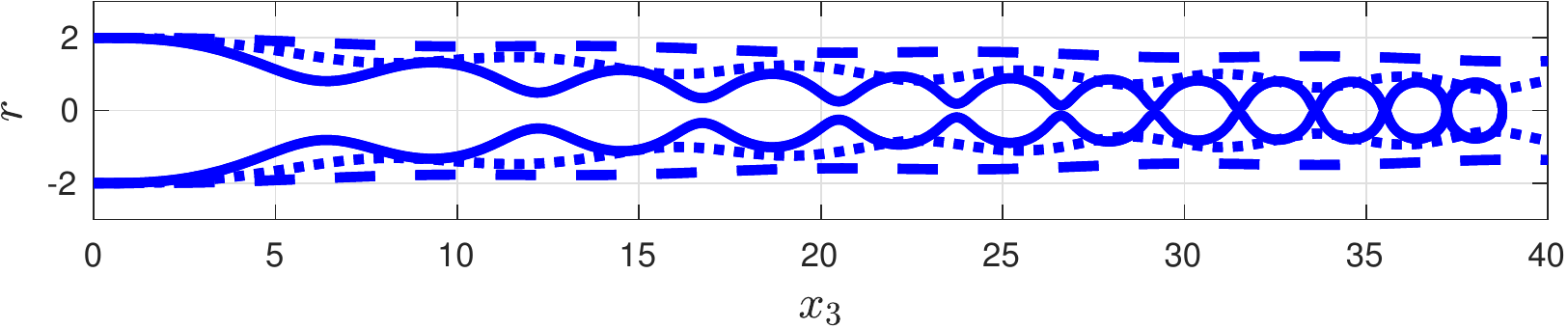}
    \caption{Cross sections of axisymmetric equilibrium pore profiles obtained from \eqref{e:IVP}, for $x_3>0$. The pores transition from slowly tapered to pinch-off for thermal gradients of $1.5^\circ$C/m (dashed), $3.2^\circ$C/m (dotted), and $15.0^\circ$C/m (solid). The simulation with the $15.0^\circ$C/m gradient terminates with a pinch-off singularity, $r(x_3)=0,$ at $x_3 =  38.7$ mm.}
    \label{fig:IVP}
\end{figure}

\subsection{Impact of Thermal Gradients on Equilibrium Pore structures}
We neglect the minor role of density stratification, setting $\delta_g=0$, so that $N_0$ is spatially constant over the inclusion, taking a value that depends only upon the volume of the inclusion and its total salt $N_T$. The steady-state version of the axisymmetric formulation shows the influence of thermal gradients on equilibrium pore shapes. Letting $r=r(x_3)$ denote the radius of an open-mouthed axisymmetric pore, then at equilibrium $\nV=0$, and combining equations \eqref{e:S-NV} and \eqref{e:curvature},  the radius satisfies the $2^{nd}$ order ODE
\beq
\label{e:IVP}
    -\frac{1+(r')^2 -r r''}{2r\sqrt{(1+(r')^2)^3}} = \talpha N_0 +\tbeta(\Theta_0(s)-\theta_*).
\eeq
The salt-free freezing temperature is $\theta_*=0^\circ$C.
The temperature profile $\Theta_0$ depends on the depth as in \eqref{e:Theta0-steady}
where $a_0$ denotes the temperature at the location $x_3=0.$ The thermal gradient $b_0$ is constant in space, but it taken at different seasonal values between $0.0015$ and $0.015$, corresponding to temperature gradients of $1.5^\circ$C/m to $15^\circ$C/m. Smaller values of $a_0$ denote a deeper sample location within the ice.

 We fix the local temperature $a_0$ and use Matlab ODE45 subroutine to resolve the
 system \eqref{e:IVP} with initial data $r'(0)=0$ and $r(0)=2$mm, corresponding to a large pore. The value of the constant salt density $N$ is selected so that $r''(0)=0$. This yields a spatially constant solution $r(x_3)=r(0)$ in the absence of a thermal gradient $(b_0=0)$. Reintroducing the thermal gradient, we solve the system for $x_3>0,$, corresponding to upwards towards the ice-air interface, the colder temperatures induce higher curvatures, and a smaller pore diameter.  Figure\,\ref{fig:IVP} shows the progression of equilibrium pore cross-sections under thermal gradients of $1.5^\circ$C/m, $3.2^\circ$C/m, and $15^\circ$C/m. The pore profile changes from weakly tapered, to oscillatory with faster tapering at intermediate gradient, to a pinch-off singularity at the largest gradient. While the oscillatory equilibrium and pinch-off states are surely unstable dynamically, they afford intuition to the role of thermal gradients in the system.

 \subsection{Impact of Thermal Gradients on Inclusion Evolution.}
 \begin{figure}
    \centering
    \includegraphics[height=2.1in]{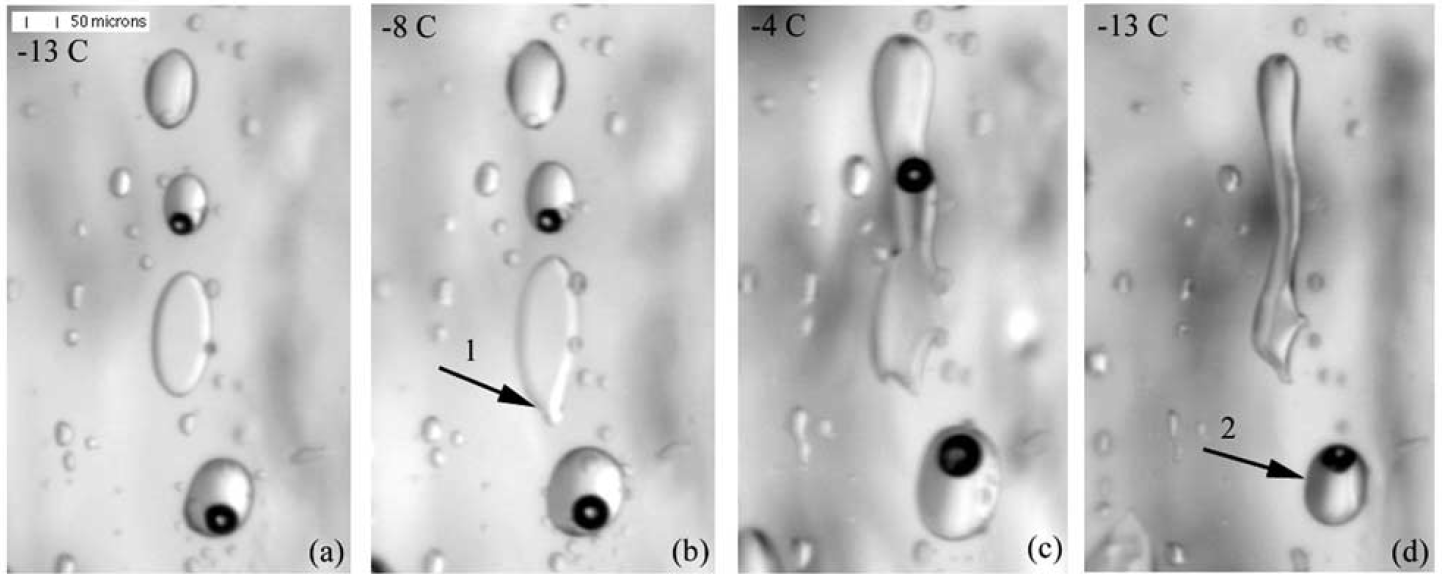}

    \caption{Evolution of brine inclusions in the absence of thermal gradient as temperature was increased from (a) $-13^\circ$C to (b) $-8^\circ$C to (c) $-4^\circ$C, then decreased to (d) $-13^\circ$C. Dark circles with bight center are air bubbles. \cite{light2003effects} reprinted with permission from John Wiley and Sons.
}
    \label{fig:light2003}
\end{figure}
\begin{figure}
    \centering
    \begin{tabular}{cccc}
    \includegraphics[width=2.6in]{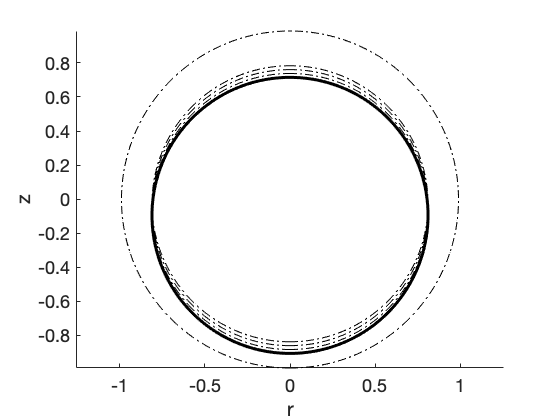} &
    \includegraphics[width=1.1in]{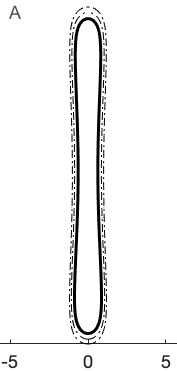} &
      \includegraphics[width=1.1in]{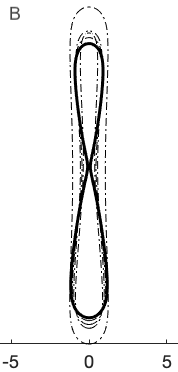} &
       \includegraphics[width=1.1in]{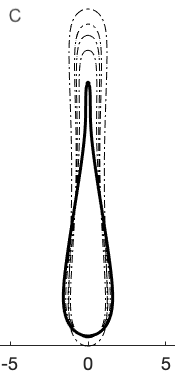}
      \end{tabular}

    \caption{Simulation of brine inclusions showing initial and intermediate times (dotted), and final time (solid). Both horizontal and vertical axis are in millimeters, but the axis in (right) are scaled compared to those in (left). Temperature and gradient are shifted at slow time $\tau=1$ in all simulations. (Left) A spherical inclusion with temperature shifted from $-2^\circ$C to $-4^\circ$C and thermal gradient shifted from $2^\circ$C/m to $14^\circ$C/m. Images are at slow time $\tau=0, 2, 3, 4, 5.$ (A) Simulation is in the absence of thermal gradient with temperature shifted from $-2^\circ$C to $-6^\circ$C. Images at $\tau=1,3, 4, 5.$ (B) Simulation has shift from 0 to $14^\circ$C/m thermal gradient and temperature decreased from $-2^\circ$C to $-10^\circ$C. Images at $\tau=0, 1, 1.5, 2, 2.58$ (pinch-off). (C) Simulation has shift from 0 to $14^\circ$C/m thermal gradient and temperature is decreased from $-2^\circ$C to $-4^\circ$C. Images at $\tau=1, 3, 4, 4.87$(pinch-off).}
    \label{fig:Wetton}
\end{figure}

The quasi-steady  system \eqref{e:asym-NV}-\eqref{e:NV-s} allows an investigation of the role of temperature gradients in the evolution of axisymmetric brine inclusions. These include both spherical  and closed cylindrical inclusions.  Spherical shapes are generically stable under curvature driven flows, indeed it is well known that motion by curvature produces spherical collapse states. In sea ice the salt concentration of a spherical inclusion of radius $R$ increases like $R^{-3}$ with decreasing radius, while curvature increases like $R^{-1}$. The build-up of salt arrests the collapse through the cryoscopic relation. As presented in Figure\,\ref{fig:Wetton} (left), we simulate a spherical brine inclusion that is initially at $-2^\circ$C with a thermal gradient of $2^\circ$C/m. At $\tau=1$ the system is subject to a rapid decrease in temperature to $-4^\circ$C and an increase in thermal gradient of $14^\circ$C/m.
In Figure\,\ref{fig:Cryoscopic}, this corresponds to an inclusion located at about $1$m depth transitioning from the June 2013 to the March 2013 temperature profiles. The inclusion contracts under the reduction in temperature, but remains largely spherical despite the asymmetry of the temperature gradient. However the gradient induces a slow rigid-body descent corresponding to a multi-dimensional traveling wave solution with a velocity that is linear in the thermal gradient for values relevant to sea ice.

Figure\,\ref{fig:light2003} presents images from an experimental investigation of reversibility of inclusion shapes under heating cooling cycles in the absence of thermal gradients. In frame (a) the unmodified first-year sea ice has a number of vertically aligned inclusions and is held at a spatially uniform $-13^\circ$C. The sample is uniformly heated with the outcome presented in frames (b) and (c), and then cooled back to $-13^\circ$C in frame (d). The four largest inclusions, initially ranging between 0.5 and 1 mm in diameter, increase in size, with one merging with a small inclusion at $-8^\circ$C, see arrow 1. At $-4^\circ$C three of the large inclusions merge into an extended brine tube of length $2$mm. Under reduction of temperature back to $-13^\circ$C, the tube contracts but is otherwise is stable.  A fundamental question is if the tube would be stable under the cooling in the presence of a thermal gradient. Indeed is is plausible that the three isolated inclusions from which the tube formed arose through the pinch-off of a tube during a cooling event in the presence of a thermal gradient.  Arrow 2 indicates an isolated pocket with reversible changes under the heating and cooling cycles.

We compare the zero-gradient experimental investigates with simulations of the quasi-steady axisymmetric Stefan-type system. While this system does not support mergings, it can follow inclusion evolution up to a pinch-off singularity. Each image in Figure\,\ref{fig:Wetton} (right - A, B, C) presents the results of three simulations, each depicting the evolution of an axisymmetric tube. In each simulation an initial inclusion was allowed to equilibrate at a given temperature under zero thermal gradient for $\tau\in(0,1)$ before being exposed to a decrease in temperature and an increase in temperature gradient.  The first simulation (A) shows the impact of a shift in temperature from $-2^\circ$C to $-6^\circ$C in the absence of a thermal gradient. This mimics the laboratory conditions of Figure\,\ref{fig:light2003}(d). The brine tube contracts but does not pinch-off. Following the cooling the evolution slows and the last two brine tube profiles are indistinguishable. The second simulation (B) shows the impact of a simultaneous application of a thermal gradient of $14^\circ$C/m and a shift in temperature from $-2^\circ$C to $-10^\circ$C. This corresponds roughly to a June to March transition, see Figure\,\ref{fig:Cryoscopic}, at a depth of $0.3$m. The brine tube pinches off quickly in the middle. This suggests that a continuation beyond pinch-off would lead to the formation of a string of two or more larger brine inclusions. Significantly there is no observable downward motion in the presence of the thermal gradient, suggesting that the spatial variation of interfacial curvature at least temporarily suppresses the transition to a traveling structure. The third simulation (C) shows the impact of  a simultaneous application of a thermal gradient of $14^\circ$C/m and a shift in temperature from $-2^\circ$C to $-4^\circ$C, corresponding to June to March transition at a depth of $0.8$m, near the bottom of the ice sheet. Here the evolution is slower, and top of the brine tube contracts into a uniformly thin neck. Significantly the pinch-off initiates at the top of the tube, which suggests that a continuation beyond pinch-off would lead to a sequence of pinch-off events that produce a string of many small brine inclusions and a single large inclusion.

\section{Discussion}

We present a thermodynamically consistent model for the slow evolution of brine inclusions within sea ice that generates salt exclusion via the entropy of salt relative to liquid water. Adapting the classical sharp-interface scaling, a multiscale analysis reduces the flow to a Stefan-type problem that couples the temperature and salinity to the evolution of the inclusion boundary. Numerical simulations of the quasi-steady version of the Stefan-type problem highlight the role of thermal gradients in the pinch-off of brine pores into spherical inclusions.
Large gradients and warm temperatures, such as found near the bottom of the ice, may lead to pore pinch-off into a range of small and large inclusions. Large gradients and swings to colder temperatures, such as found at the top of the sea ice, may lead to pinch-off into equal sized inclusions.  These observations are commensurate with  Figure\,\ref{fig:3D-brine} which displays samples of brine inclusions in first-year sea ice. In both images the brine inclusions arise in vertically aligned columns that are evocative of the pinch-off of longer brine pores examined here.
It is also natural to ask if  thermal gradients will induce migration of brine inclusions. Interestingly, the analysis and numerics suggest that spherical inclusions are more susceptible to migration, as their constant curvature does not readily adjust to the inhomogeneity of the thermal gradient, leading to rigid-body evolution rather than the deformation that seen in pore-type inclusions.

The model presents several opportunities for novel analysis. The combination of an $L^2$ gradient for the non-conserved phase, a weighted $H^{-1}$ gradient for the conserved salt leads, and chemotactic terms arising from the relative entropy present several challenges to the analysis. Most chemotaxis results have addressed spatially localized patterns, such as spikes, \cite{WW05}. Brine inclusions are very much a chemotaxis phenomena, transporting uniformly distributed ocean salt at 3.5\% weight fraction into brine inclusions at 10-20\% weight fraction. However the inclusions are spatially extended patterns, with a length scale that is $10^6$ times longer than the ice-liquid interface. A rigorous analysis of the stability and evolution of these spatially extended inclusions seems to require a different class of tools. A good starting point is to address the stability of traveling spherical inclusions in the presence of a thermal gradient.

The issues raised present opportunities for model calibration. Almost all laboratory work on sea ice is conducted at constant temperature, without thermal gradient. The thermal gradients in sea ice in winter are very significant. A simple point of validation would be to measure the drift speed of circular inclusions as a function of the strength of the thermal gradient, or to recapitulate the work of \cite{light2003effects}, such as that presented in Figure\,\ref{fig:light2003},  under the influence of thermal gradients.

The model contains many simplifications, some of which make it harder to incorporate experimental data into the initial model development.  An obvious improvement is to consider a more physical balance between entropy and free energy for the ice-liquid transition. As discussed in the work of Penrose and Fife, \cite{penrose1990thermodynamically}, the entropy of the phase change should be convex, with the non-convexity that drives the spinodal decomposition appearing through the temperature dependence of the latent heat. An even more ambitious extension is to incorporate the microstructure of the ice phase and the elastic energy driven by the expansion of water upon freezing. This would require the full GENERIC framework, \cite{mielke2011formulation}.
Both of these projects are future work.

\section{Acknowledgment}
KP recognizes the support of the National Science Foundation through grant DMS 1813203. BW acknowledges support from an NSERC Canada grant.

\bibliographystyle{abbrv}
\bibliography{BrineInclusionsBib}
\end{document}